\documentclass[12pt,a4paper]{article}
\usepackage{euscript,amsfonts,amssymb,amsmath,amscd}
\newcommand{\op}{\operatorname}
\begin{document}
\begin{center}
{\bf Cobordism realizations of Thom polynomials of singularities.}\\

Andrei Kustarev.
\end{center}

Let $\EuScript E^k(m,n)$ be a vector space of k-jets of holomorphic maps $(\mathbb C^m,0)\to(\mathbb C^n,0)$. The group $A(m,n):=Hol(\mathbb C^m,0)\times Hol(\mathbb C^n,0)$ acts on $\EuScript E^k(m,n)$ by $(\phi,\psi)\cdot f=\psi\circ f\circ\phi^{-1}$. By singularity $\eta$ we understand an irreducible $A(m,n)$-invariant complex subvariety in $\EuScript E^k(m,n)$. A holomorphic map $f:M^{2m}\to N^{2n}$ has the singularity $\eta$ in point $x\in M$, if the k-jet of $f$ in $x$ belongs to $\eta$ (after writing it in any local coordinates).

The following result is well-known in singularity theory. Let $\eta_f$ be a set of points, where $f$ -- a general holomorphic map -- has the singularity $\eta$. Then cohomology class $[\eta_f]$ can be expressed as a polynomial $P_H^{\eta}(c_i(TM),c_j(f^*TN))$ depending on cohomological characteristic classes $c_i(TM),c_j(f^*TN)$.  $P_H^{\eta}(c_i(TM),c_j(f^*TN))$ is called Thom polynomial, or universal polynomial of singularity $\eta$, because its coefficients depend only on certain type of $\eta$ and not on $f$.

{\bf Definition 1.} A polynomial $P_U^{\eta}(a_i,b_j)\in\Lambda[a_i,b_j]$ (here $\Lambda=U^*(pt)$) is called a cobordism realization of Thom polynomial $P_H^{\eta}(a_i,b_j)$, if $\epsilon(P_U^{\eta}(c_i^U(TM),c_j^U(f^*TN)))$=$P_H^{\eta}(c_i(TM),c_j(f^*TN))$ for any $f$. Here $\epsilon:U^*(X)\to H^*(X)$ is a canonical map, $c_i^U(\xi)$ is the i-th Chern class of $\xi$ in cobordisms.

An obvious realization is $P_U^{\eta}(a_i,b_j)\equiv P_H^{\eta}(a_i,b_j)$. The consequence is that Steenrod problem of realization of cycle $\eta_f$ by image of a smooth manifold always has a solution.

The cycle $\eta_f$ is a complex semialgebraic subset in $M$. Denote by $\tilde\eta_f$ the set of its regular points. 

{\bf Definition 2.} A realization $P_U^{\eta}$ will be called nice if the following conditions hold:
\begin{enumerate}
{\item (Resolution). The class $P_U^{\eta}(c_i^U(TM),c_j^U(f^*TN))$ can be represented by complex-oriented map $g:X\to M$ such that $\op{Im} g=\eta_f$ and $\tilde\eta_f$ are regular values of $g$.}
{\item (Functoriality). Suppose that $h:N_1\to N$ is a holomorphic map transversal to $M$. Then there exists a complex manifold $M_1$ and holomorphic maps $f_1:M_1\to N_1$ and $\tilde h:M_1\to M$ such that $h\circ f_1=f\circ\tilde h$. 
We require $\tilde h^*P_U^{\eta}(c_i^U(TM),c_j^U(f^*TN))$=$P_U^{\eta}(c_i^U(TM_1),c_j^U(f_1^*TN_1))$.}
\end{enumerate}

{\bf Disappointing Lemma.} The realization $P_U^{\eta}(a_i,b_j)\equiv P_H^{\eta}(a_i,b_j)$ isn't nice.

Consider a general algebraic map $f:CP^2\to CP^2$ of degree $d^2, d>1$. Let $X\subset CP^2$ be a set of points where $\dim\ker df>0$ (the $\Sigma^1$ singularity). Denote by $\nu$ the dual to canonical linear bundle over $CP^2$. Then $X$ is a non-singular curve and $[X]=c_1^U(\nu^{\otimes(3d-3)})$ in $U^2(CP^2)$.

The cohomological Thom polynomial of $\Sigma^1$ for general maps $f:M^{2n}\to N^{2n}$ is $c_1(f^*TN)-c_1(TM)$ (see \cite{thom},\cite{port}). If $c_1^U(f^*TN)-c_1^U(TM)$ is a nice realization, then one can show that $[X]=(3d-3)c_1^U(\nu)$ in $U^2(CP^2)$, using condition 1. This contradiction proves lemma.

The set of nice realizations is not empty. For every $\eta$, one can construct at least one nice realization using Hironaka resolution theorem and universal vector bundle with a fiber $\EuScript E^k(m,n)$. 

Even in the simplest case of $\Sigma^1$ singularity of maps $f:M^{2n}\to N^{2n}$ this set consists at least of two elements. One nice realization is $P_1(\Sigma_1):=c_1^U(\det(f^*TN-TM))$, where $\det\xi$ is the linear bundle with $c_1(\det\xi)$=$c_1(\xi)$.

The second nice realization $P_2(\Sigma^1)$ can be constructed by using general resolution of $\Sigma^r$ (see \cite{kaz}). One can prove that $P_1(\Sigma^1)$ and $P_2(\Sigma^1)$ are different by considering a general map $CP^4\to CP^4$ and comparing the cobordism classes of $P_1(\Sigma^1)$ and $P_2(\Sigma^1)$.

If $P_U^{\eta}$ and $Q_U^{\eta}$ are two nice realizations of $\eta$, then any integer combination of the form $\lambda P_U^{\eta}+(1-\lambda)Q_U^{\eta}$ is also nice realization. Denote by $\bar\eta$ the set of singular points of $\eta\subset\EuScript E^k(m,n)$ and by $P_H^{\bar\eta}(a_i,b_j)$ its Thom polynomial. 

{\bf Theorem.} Let $P_U^{\eta}(a_i,b_j)$ and $Q_U^{\eta}(a_i,b_j)$ be two nice realizations of $\eta$. Then $\op{ch}_U (P_U^{\eta}-Q_U^{\eta})$ is divisible by $P_H^{\bar\eta}(a_i,b_j)$. Here $\op{ch}_U$ is the Dold-Chern character in complex cobordisms (see \cite{buch}).

The last two statements may be useful in classification problem of nice realizations.

The author is grateful to Prof.V.M.Buchstaber for constant attention to this work and to Prof.M.E.Kazarian, G.A.Merzon and E.A.Gorsky for useful discussions.

\end{document}